\documentclass{amsart}
\numberwithin{equation}{section}
\usepackage{amssymb}
\let\blb\mathbb

\def\QQ{{\blb Q}}

\def \ZZ{{\blb Z}}

\def\p{\mathfrak{p}}
\def\m{\mathfrak{m}}

\def\Lotimes{\overset{L}{\otimes}}

\def\Mod{\operatorname{Mod}}
\def\Char{\operatorname{char}}

\def\rad{\operatorname {rad}}

\def\Spec{\operatorname {Spec}}

\def\Ext{\operatorname {Ext}}
\def\Hom{\operatorname {Hom}}

\def\End{\operatorname {End}}
\def\RHom{\operatorname {RHom}}

\def\GKdim{\operatorname {GKdim}}

\def\End{\operatorname {End}}

\def\rk{\operatorname {rk}}

\def\gldim{\operatorname{gl\,dim}}

\def\r{\rightarrow}

\DeclareMathOperator{\Div}{Div}

\theoremstyle{plain}
\newtheorem{lemma}[equation]{Lemma}
\newtheorem{proposition}[equation]{Proposition}
\newtheorem{theorem}[equation]{Theorem}
\newtheorem{corollary}[equation]{Corollary}

\theoremstyle{definition}

\newtheorem{question}[equation]{Question}

{

}

\newtheorem{remark}[equation]{Remark}
\newtheorem{remarks}[equation]{Remarks}

\newdimen\uboxsep \uboxsep=1ex
\def\uboxn#1{\vtop to 0pt{\hrule height 0pt depth 0pt\vskip\uboxsep
\hbox to 0pt{\hss #1\hss}\vss}}
\def\uboxs#1{\vbox to 0pt{\vss\hbox to 0pt{\hss #1\hss}
\vskip\uboxsep\hrule height 0pt depth 0pt}}

\begin{document}
 \title{Noncommutative resolutions and rational singularities}
\author{J. T. Stafford}
\address{Department of Mathematics, East Hall, 530 Church Street,
University of Michigan, Ann Arbor,
MI 48109-1043, USA.} 
\email{jts@umich.edu}
\author{M. Van den Bergh}
 \address{Departement WNI,  Universiteit Hasselt,
 3590 Diepenbeek, Belgium.}
  \email{michel.vandenbergh@uhasselt.be}
  \thanks{The  first author was
  partially supported by the NSF through grants DMS-0245320 and 
  DMS-0555750  and also  by the  Leverhulme Research Interchange Grant
F/00158/X. Part of this work was written up while the first
 author was visiting and supported by the Newton Institute, Cambridge. He would 
 like to thank all three institutions for their financial support.}
\thanks{The second author is a senior researcher at the FWO}
\keywords{noncommutative
 geometry, homologically homogeneous algebras, rational singularities}
\subjclass{14A22, 14E15, 16S38, 18G20}
\begin{abstract} Let $k$ be an algebraically closed field of
  characteristic zero. We show that the centre of a homologically
  homogeneous, finitely generated $k$-algebra has rational
  singularities. In particular if a  finitely generated normal commutative 
  $k$-algebra has a noncommutative crepant resolution, as introduced by the 
  second author,  then it has rational singularities.
\end{abstract}

\maketitle


\section{Introduction}\label{ref-1-0}

Throughout the paper, $k$ will denote a fixed algebraically closed 
field of characteristic zero  and unless otherwise  specified
all rings will be $k$-algebras.  Suppose that $X=\Spec R$ for 
an affine (that is, finitely generated)  normal Gorenstein  $k$-algebra $R$.
  The nicest form of resolution of singularities  
 $f:Y\to X$ occurs when $f$ is  \emph{crepant} in the sense that 
$f^\ast \omega_X=\omega_Y$.  Even when they exist, crepant resolutions 
 need not be unique, but they are related---indeed Bondal and Orlov conjectured 
 in \cite{BO1} (see also  \cite{BO2}) that two such resolutions should be derived equivalent. 

  Bridgeland \cite{Bridge} proved the Bondal-Orlov
 conjecture in dimension~$3$. The second author observed in
 \cite{VdB-crepant} that Bridgeland's proof could be explained in
 terms of a \emph{third} crepant resolution of $X$ which is now
 noncommutative (the definition will be given below) and this had lead to a number of different
 approaches to the Bondal-Orlov conjecture and related topics---see, for example,
  \cite{Bezrukavnikov, Bezkal2, Chen,  IR,  Kaledin, Kawamata1}.

  It is therefore natural to ask how the existence of a noncommutative crepant resolution 
  affects the original commutative singularity.
 It is  well-known,  and follows easily from \cite[Theorem~5.10]{KM},   that if
 a Gorenstein singularity has a crepant resolution then it has
 rational singularities. So it is logical to ask, as was raised in  \cite[Question~3.2]{VdB32},
  {\em  is this  true for  a noncommutative
   crepant resolution}? In this paper we answer this question
 affirmatively, but  before stating the result precisely, we need to define the relevant terms.

\medskip  
Let $\Delta$ be a prime  affine $k$-algebra that is finitely generated  as a module 
over its centre $Z(\Delta)$.   Mimicking \cite{BH}, we say that $\Delta$ is \emph{homologically
  homogeneous of dimension $d$} if all simple $\Delta$-modules have
the same projective dimension $d$.  By \cite{Ra} and \cite{BH} such a
 ring  $\Delta$ has global and Krull dimensions  equal to $d$ and, as has   been
shown in \cite{BH}, the properties of homologically homogeneous rings
closely resemble those of commutative regular rings. So the idea is to use such a ring 
$\Delta$ as a noncommutative analogue of a crepant resolution. Formally, following 
\cite{VdB32} we define a \emph{noncommutative crepant resolution}  of
$R$ to be any homologically homogeneous ring of the form
$\Delta=\End_R(M)$, where $M$ is a reflexive, finitely generated
$R$-module. We refer the reader to \cite[Section~4]{VdB32} for the logic behind this definition.

 The main result of  the present note is  the following: 

\begin{theorem} \label{ref-1.1-1} {\rm (Theorem~\ref{ref-4.1-31})}
 Let $\Delta$ be a  homologically homogeneous  $k$-algebra.  Then  the centre 
  $Z(\Delta)$ has rational singularities.

In particular if a normal affine $k$-domain  $R$  has a noncommutative 
crepant resolution then it  has rational singularities.
\end{theorem}

\medskip

In Section~\ref{ref-5-34}
we give two examples related to the theorem. The first example shows
that if $\Delta=\End_R(M)$ has finite global dimension  then it need not be homologically
homogeneous even under reasonable hypotheses on $M$ and $R$.
The second shows that  Theorem~\ref{ref-1.1-1} can fail in positive characteristic. 

\medskip 

\noindent
{\bf Notation.}
Throughout the paper $R$ will be a normal commutative noetherian $k$-domain 
and  $\Delta$   will be a $k$-algebra, with 
centre $Z=Z(\Delta)$ containing $R$, 
  such that   $\Delta$ is a finitely generated $R$-module. We say that 
  $R$ is \emph{essentially affine} if it is a localization of an affine $k$-algebra.
The dimension function used in this paper will be the \emph{Gelfand-Kirillov dimension} 
of $\Delta$ as a $k$-algebra, written $\GKdim \Delta$.
 By \cite[Proposition~8.2.9(ii) and Theorem~8.2.14(ii)]{MR} 
$\GKdim \Delta = \GKdim R$ and $\GKdim R$ 
is just the transcendence degree of $R$ over $k$.


\section{Homologically homogeneous rings}\label{ref-2-2}

In this section we introduce homologically homogeneous rings 
and prove some basic facts about their structure and their dualizing 
complexes.  Many of these results use the machinery of tame orders 
and so we    start by discussing this concept.

\subsection*{Tame orders}  \label{ref-2-3}   
Assume that  $\Delta$  is \emph{a   prime $R$-order in $A$}, by which we mean that 
$\Delta$ is a prime ring with simple artinian ring of fractions $A$. 
We write $\mathfrak{P}_1=\mathfrak{P}_1(R) $
for the set of height one prime ideals of $R$ and say that \emph{a
  property $\mathcal{P}$ holds for $\Delta$ in codimension one} if it
holds for all $\Delta_\p=\Delta\otimes_RR_\p :
\p\in\mathfrak{P}_1$.  Following \cite{silver}, the prime $R$-order 
$\Delta$ is called a
\emph{tame $R$-order} if $\Delta$ is a finitely generated, reflexive
$R$-module that is   hereditary  in codimension one.

 The paper \cite{silver} 
implicitly assumes that $R=Z(\Delta)$, but we prefer  not make this assumption.
However, by the following standard result,   the question of whether $\Delta$ is 
a tame $R$-order is independent of the choice of normal central subring  $R$.

 \begin{lemma}\label{ref-2.2-5} Let $\Delta$ be a tame $R$-order. Then a finitely generated
 $\Delta$-module is reflexive as an $R$-module if and only if it is reflexive as a $\Delta$-module. 
\end{lemma}

\begin{proof}  
By \cite[Corollary~1.6]{silver} (which does not require $R=Z(\Delta)$)
a $\Delta$-reflexive module is $R$-reflexive. Conversely, 
suppose that $M$ is a finitely generated $\Delta$-module that is $R$-reflexive.
Since $M$ is 
therefore torsion-free as a $\Delta$-module,  $M_\p=M\otimes_RR_\p$ is torsion-free and hence 
projective   over the hereditary prime ring $\Delta_\p$,  for all $\p\in\mathfrak{P}_1$.   Thus, 
by \cite[Lemma~1.1]{silver}, 
$$\begin{array}{rl}
M\ = & \bigcap_{\p\in \mathfrak{P}_1} M_\p = 
 \bigcap_{\p\in \mathfrak{P}_1}\Hom_{\Delta_p}(\Hom_{\Delta_p}(M_\p,\,\Delta_\p),\,\Delta_\p)\\
  \noalign{\vskip 5pt}  &   \supseteq
  \Hom_{\Delta}(\Hom_{\Delta}(M,\,\Delta),\,\Delta).
 \end{array}$$
Thus, $M= \Hom_{\Delta}(\Hom_{\Delta}(M,\,\Delta),\,\Delta)$, as required.
\end{proof}

Let $\Delta$ be a tame $R$-order in $A$. A \emph{divisorial}
fractional $\Delta$-ideal is  any  reflexive fractional $\Delta$-ideal in $A$ 
 that is invertible in codimension one. 
  By \cite[Theorem~2.3]{silver}, divisorial
fractional ideals form a free abelian group $\Div(\Delta)$ with
product $I\cdot J=(IJ)^{\ast\ast}$ where $K^\ast=\Hom_R(K,R)$ denotes the
$R$-dual of a fractional ideal $K$.  The $n^{\mathrm{th}}$  power  
 $(I^n)^{\ast\ast} $  of $I$ under this dot-operation  is called  
 the \emph{$n^{\mathrm{th}}$  symbolic power of $I$} and written $I^{(n)}$.
Write $\rad S$ for the Jacobson radical of a ring $S$. 


\subsection*{Homologically homogeneous rings}
Homologically homogeneous ring,  as defined in the   introduction, 
 have a particularly pleasant structure and the following result 
 provides some of the  properties we will need. 

\begin{theorem} 
\label{ref-2.3-6}
Assume that $\Delta$ is   homologically homogeneous of dimension $d$. 
\begin{enumerate}
\item $\Delta$ is    CM as a   module over its centre $Z$.
\item 
Both $\GKdim \Delta$  and the global homological dimension $\gldim \Delta$ of $\Delta$ 
equal~$d$.
\item $Z$  is an affine  CM   normal   domain. 
\item $\Delta$ is a   tame $Z$-order.
\end{enumerate}
\end{theorem}

\begin{proof}   (1,2) By \cite[Theorem~8]{Ra}, $\gldim \Delta= d$.
The rest of  parts (1) and (2) follow from \cite[Theorem~2.5]{BH}.

  (3)  By hypothesis, $\Delta$ is finitely generated as both a $Z$-module and a $k$-algebra. Thus
 the Artin-Tate Lemma \cite[Lemma~13.9.10]{MR} implies that  $Z$ is an affine $k$-algebra.
As $\Char k=0$, the reduced trace map $\Delta\to Z$ is surjective and 
  so $Z$ is a $Z$-module summand of $\Delta$. Thus  $Z$ is CM by part (1).
  As $\Delta$ is prime, $Z$ is a domain, while $Z$ is normal by  
\cite[Theorem~6.1]{BH}.
 
(4)  As  $\Delta$ is CM as a $Z$-module, it  is certainly reflexive.
  By \cite[Corollary~2.2 and Theorem~2.5]{BH},   $\Delta$ is hereditary in  codimension one.  
\end{proof}

The standing  assumption that $k$ have  characteristic zero is crucial for the 
proof of  part (3) of the theorem.  Indeed,  \cite[Example~7.3]{BHM} shows that
the centre  $Z(\Gamma)$ of a homologically homogeneous ring $\Gamma$ 
 need not be CM in bad characteristic.
 
The following criterion for a ring to be homologically homogeneous will be useful.

\begin{lemma}  \label{ref-2.4-7}  
Suppose  that  $R$  is  an affine $k$-algebra and that   $\Delta$ is a
prime ring. If  $\Delta$ is  a CM $R$-module  with
  $\GKdim \Delta =\gldim \Delta$, then $\Delta$ is homologically homogeneous.
\end{lemma}

\begin{proof}  This is, essentially, \cite[Proposition~7.2]{BH}, but here is a direct proof.
Suppose that  $S$ is a simple $\Delta$-module with projective dimension 
  $u <d=\gldim \Delta$ and consider a projective   $\Delta$-resolution of $S$:
\[
0\r P_u\r \cdots \r P_1\r P_0\r S\r 0.
\]
Viewed as a complex over $R$ this is a  resolution of length $<d$ of a finite
length $R$-module by CM modules of dimension $d$. 
An  easy  depth argument shows that this is impossible. 
\end{proof}


\subsection*{Dualizing modules and complexes} 
 In order relate properties of  a homologically homogeneous ring to those of its centre
  we will use the machinery of dualizing complexes and  we 
discuss their structure   in this subsection. 
Most of the background material   comes from \cite{VdB16, ye5, YZ1, YZ2}  
and the reader is referred  to those papers for more details.   
 Throughout this discussion, in addition to our standing assumptions, we assume
  that $R$ is essentially affine.

 Write  $\Delta^e=\Delta\otimes_k\Delta^{\mathrm{op}}$ and denote the  derived category
 of left $\Delta^e$-modules by $D(\Delta^e) $.  Following
 \cite{Ye}, a \emph{dualizing complex} for $\Delta $ is a complex of
 $\Delta$-bimodules $D$, with finite injective dimension on both
 sides, such that
  \begin{enumerate}
  \item the cohomology of $D$ is given by $\Delta$-bimodules that are finitely 
  generated on both sides, and
  \item   in   $D(\Delta^e)$  
  the pair of natural morphisms  $\Phi: \Delta\to \mathrm{RHom}_\Delta(D,D)$ and 
  $\Phi^o : \Delta \to  \mathrm{RHom}_{\Delta^{\mathrm{op}}}(D,D)$ are isomorphisms. 
    \end{enumerate}
Following \cite[Definition~8.1]{VdB16}, the 
dualizing complex $D_\Delta$ is called \emph{rigid} if there is an isomorphism 
$\chi: D_\Delta\cong \mathrm{RHom}_{\Delta^e}(\Delta, \,D_\Delta\otimes  D_\Delta)$
in $D(\Delta^e)$.  
The significance of rigidity is that, although dualizing complexes are not unique, 
rigid dualizing complexes are, in the sense that the pair $(D_\Delta,\chi)$ is unique up to a 
unique isomorphism   \cite[Proposition~8.2]{VdB16} \cite[Theorem~3.2]{YZ1}.

Although dualizing complexes (rigid or otherwise) do not exist for all
finitely generated  noncommutative noetherian rings \cite[p.~529]{KRS},  by 
 \cite[Proposition~5.7]{ye5} and \cite[Theorem~3.8]{YZ2}  they do exist for our 
 rings $R$ and  $\Delta$. 

Write $d=\GKdim \Delta=\GKdim R$.
 The cohomology of $D_R$ and $D_\Delta$ lies in degrees $\ge -d$ and
 we define $\omega_R=H^{-d}(D_R)$ and
 $\omega_\Delta=H^{-d}(D_\Delta)$.  An important fact
 \cite[Corollary~3.6]{YZ1} is that the cohomology of $D_\Delta$ is
 \emph{$Z$-central} in the sense that the left and right actions of
 $Z$  agree.  In particular, $\omega_\Delta$ is $Z$-central.
 
  The following results gives some basic properties that we will need about these objects.
  If $M$ is $\Delta$-bimodule then $Z(M)=\{w\in M : \delta w=w\delta 
\text{ for all } \delta \in \Delta\}$ is called   \emph{the centre    of~$M$.} 

  \begin{lemma} \label{new-lemma} 
Assume that   $R$ is an  essentially affine $k$-algebra.  Then:
  \begin{enumerate}
    \item\label{new-lemma1}  $D_{\Delta}\cong\RHom_{R}(\Delta,D_R)$ in
 $D(\Delta^e)$.
    \item \label{new-lemma2}
    $\omega_\Delta\cong\Hom_R(\Delta,\omega_R)$ as $\Delta^e$-modules.
    \item\label{new-lemma3}
      If $\mathcal{C}\subset Z$ is  multiplicatively closed  then 
    $\omega_{\Delta_\mathcal{C}}\cong (\omega_\Delta)_\mathcal{C}$ as $\Delta$-bimodules. 
    \end{enumerate}
    Assume  in addition that     $\Delta$ is a tame $R$-order. Then: 
         \begin{enumerate}
     \item[(4)] \label{new-lemma4}
     $\omega_\Delta$ is a reflexive   as a left or right  $\Delta$-module.
    \item[(5)] \label{new-lemma5}
 $\omega_\Delta$ is invertible in codimension one. Moreover, as bimodules
$$   
\omega_\Delta=\biggl( \omega_{Z} \otimes_{_{\scriptstyle Z}} 
 \prod_{\p\in \mathfrak{P}_1(Z)}
 \bigl(\Delta\cap \rad(\Delta_\p)\bigr)  \cdot\p^{(-1)}\biggr)^{\ast\ast}.
 $$
  \item[(6)]\label{new-lemma6}   There is a canonical isomorphism
$Z(\omega_\Delta)=\omega_Z$.
    \end{enumerate}
    \end{lemma}

    \begin{proof}  
    (1) The proof of \cite[Proposition~5.7]{ye5} shows that $\RHom_{R}(\Delta,D_R)$ is a rigid 
    dualizing complex for $\Delta$ and so the result follows   by the uniqueness of $D_\Delta$. 

   (2) Take cohomology  of \eqref{new-lemma1}.

     (3)   By \cite[Theorem~3.8]{YZ2}  
 $D_{\Delta_\mathcal{C}}\cong\Delta_\mathcal{C}\Lotimes_\Delta 
 D_\Delta\Lotimes \Delta_\mathcal{C}$ as   $\Delta$-bimodules.
  Now take cohomology, using the fact that, as mentioned above, each
  $\mathrm{H}^q(D_\Delta)$  is $Z$-central. 

(4)  By part (2) and \cite[Lemma~1.5]{silver}  it suffices to prove the
 result for $\omega_R$.
This case   is well-known, but here is an easy  proof.
By part (3) we may assume that $R$ is  an affine  $k$-algebra. By
 Noether normalization $R$ is a finitely generated module over some polynomial subring
  $R_0$ and it is a tame $R_0$-order since it is normal. 
It is standard  \cite[Example~ 3.13]{YZ2} that $\omega_{R_0}
\cong R_0$  as bimodules and so 
\cite[Lemma~1.5]{silver} and Lemma~\ref{ref-2.2-5}  imply that
 $\omega_R \cong\Hom_{R_0}(R,\omega_{R_0})$ is a reflexive $R$-module. 

(5) The  first assertion follows, for example, from 
\cite[Corollary~37.9]{CR} combined with part~\eqref{new-lemma2}.  
In codimension one  the displayed equation also   follows easily from part~\eqref{new-lemma2} and 
 the general case then follows from  parts (4) and~\eqref{new-lemma3}.

(6) This follows  from part (5).
  \end{proof}

\begin{remark}\label{canonical-defn}
We emphasize that our definition of $\omega_R$ does coincide with the usual commutative
notion   $\bigwedge^d(\Omega_{R/k})^{**}$
when $R$ is  essentially affine with $\GKdim R =  d$. To see this, set  
$\omega_R'=\bigwedge^d(\Omega_{R/k})^{**}$. Then \cite[Lemma~5.4]{ye5}
shows that $\omega_S=\omega_S'$   holds for any  regular, essentially finite domain $S$.
As $R$ is normal, it is regular in codimension one and so Lemma~\ref{new-lemma}(3)
implies that $(\omega_R)_\p=(\omega'_R)_\p$ for all height one prime ideals $\p$.
By  Lemma~\ref{new-lemma}(4), $\omega_R$ 
and $\omega_R'$ are reflexive, and hence  $\omega_R=\omega'_R$.
\end{remark}

\begin{proposition} 
\label{ref-2.15-18} Assume that $\Delta$ is a  prime affine $k$-algebra. 
Then $\Delta$ is homologically homogeneous of 
dimension  $d$ if and only if $\gldim \Delta<\infty$ and
  $D_\Delta=\Omega[d]$ for some invertible $\Delta$-bimodule~$\Omega$.
 If this holds then  $\Omega=\omega_\Delta$. 
\end{proposition}

\begin{remark}\label{ref-2.16-19} 
  In the notation of \cite[Section~8]{VdB16}, the proposition states
  that $\Delta$ is homologically homogeneous of dim\-ension $d$ if and
  only if $\gldim \Delta<\infty$ and $\Delta$ is AS-Gorenstein.
See \cite[Theorems~1.3 and~1.4]{SZ} for a closely related result.
    \end{remark}

  \begin{proof} 
Assume first that $\Delta$ is homologically homogeneous of dimension $d$.  
Since the statement of the proposition is independent of the choice of $R$ 
we may, by Noether normalization, assume 
  that $R$ is a polynomial ring. By   Theorem~\ref{ref-2.3-6}(1,3), $\Delta$ is   CM and hence free 
 as an $R$-module. 
  But now $D_R=\omega_R[d]\cong R[d]$ and so   
  $D_\Delta=\RHom_R(\Delta,D_R)$ lives purely in dimension $-d$, whence  
  $D_\Delta=\omega_\Delta[d]$.
    Lemma~\ref{new-lemma}\eqref{new-lemma2} implies that $\omega_\Delta$ is free 
and hence CM as an $R$-module and so      \cite[Corollary~3.1]{BH}
implies that  $\omega_\Delta$  is a projective   $\Delta$-module on either side. On the other hand,
  as $\Delta$ is  a  free $R$-module it  is a tame $R$-order and so Lemma~\ref{new-lemma}(5)
     implies that $\omega_\Delta$ is  invertible in codimension one. Together
     with \cite[Proposition~3.1]{silver},  
      these observations imply that   $\omega_\Delta$ is invertible,  finishing the proof in this  direction. 

   Conversely, assume that $\gldim \Delta<\infty$ and
  $D_\Delta=\Omega[d]$ for some invertible bimodule 
   $\Omega$. We will show that
every simple $\Delta$-module  $S$ has projective dimension~$d$.
By \cite[Corollary~6.9]{YZ1} $D_\Delta$ is Auslander and GKdim-Macaulay
  in the sense  of \cite[Definitions~2.1 and 2.24]{YZ1}. 
 Since $S$ is finite dimensional the Macaulay property means that 
   $\Ext^d_\Delta(S,\Omega)\neq 0$.
If $\gldim\Delta =e>d$ then, by \cite[Theorem~8]{Ra},
there exists a simple $\Delta$-module $S$ with $\Ext^e_\Delta(S,\Delta)\neq 0$.
Since $\Omega$ is invertible, this implies that $E=\Ext^e_\Delta(S',\Omega)\neq 0$
for  $S'=\Omega\otimes_\Delta S$.
By the Auslander property,  this means that 
  $\operatorname{Cdim}_{D_\Delta}(E)\leq -e$  which, 
by the GKdim-Macaulay property,    implies that 
$\GKdim S<0$. This is absurd.
\end{proof}

The following formul\ae\ will be useful.

\begin{corollary}\label{ref-2.17-20}
 Assume that $R$ is essentially affine with $\GKdim R=d$ 
  and let $\Delta$  be a  tame $R$-order.  Then 
\begin{equation}
\label{ref-2.18-21}
\omega_\Delta^{(-1)}=\Ext^d_{\Delta^e}(\Delta,\Delta^e)^{\ast\ast}.
\end{equation}
If $\Delta$ is homologically homogeneous then 
\begin{equation}
\label{ref-2.19-22}
\omega_\Delta^{(-1)}=\RHom_{\Delta^e}(\Delta,\Delta^e)[d].
\end{equation}
\end{corollary}

\begin{proof}  If $\Delta$ is homologically homogeneous  then it has dimension $d$ by  
Theorem~\ref{ref-2.3-6}. Thus    \cite[Proposition~8.4]{VdB16}
and Remark~\ref{ref-2.16-19} combine to prove 
  \eqref{ref-2.19-22}.

Now suppose that   $\Delta$ is a tame $R$-order and set $\Gamma=\Delta_\p$, 
for some  $ \p\in\mathfrak{P}_1(Z)$. Then $\Gamma$ is an hereditary order and,
 by  \cite[Theorem~13.10.1]{MR},
   $Z(\Gamma)$ is a  local PID. By Lemma~\ref{new-lemma}(4,5), $\omega_\Delta$ is
 invertible and hence,   just as in the proof of Theorem~\ref{ref-2.15-18}, 
 $D_\Delta=\omega_\Delta[d]$. Thus, 
  \cite[Proposition~8.4]{VdB16} can again   be applied to show that 
  $\omega_\Gamma^{(-1)}=\Ext^d_{\Gamma^e}(\Gamma,\Gamma^e)^{\ast\ast}$.
In other words,  \eqref{ref-2.18-21} holds in codimension one.
Since both sides
of that equation  are reflexive, it  holds
everywhere.  
\end{proof}


\section{Reduction to the Calabi-Yau case} 
Let $\Delta$ be a homologically homogeneous ring.   In Section~\ref{ref-4-30} we  will use the
structure of $\omega_\Delta$ to show that $Z$ has rational
singularities, but this is awkward to prove when $\omega_\Delta$ is not
cyclic. In this section we show how to use a trick from
 \cite[Theorem~3.1]{NV} to  (locally) replace $\Delta$ by an
order for which  $\omega_\Delta$ is
generated by a single central element. This is  a noncommutative 
generalization of a well known 
technique in algebraic geometry  where one constructs a
Gorenstein cover of a $\QQ$-Gorenstein singularity. 

Given a tame $R$-order $\Gamma$ 
in $A$ and  $I\in \Div(\Gamma)$, the  \emph{Rees ring $\Gamma[I]$}
of $\Gamma$ is defined to be  the subring  $\sum_{n=-\infty}^\infty I^{(n)}x^n$ of the Laurent 
polynomial ring $ A[x,x^{-1}].$

\begin{proposition} \label{ref-3.1-23} Assume  that 
$\Delta$ is   homologically homogeneous.
 For some $n\geq 1$, suppose  that $\omega^{\otimes
    n}_\Delta\cong \Delta$ as bimodules and choose $n$  minimal with
  this property. Write 
$$
\Lambda = \Delta\oplus \omega_\Delta\oplus \omega^{\otimes 2}_\Delta\oplus\cdots\oplus
\omega^{\otimes n-1}_\Delta
$$
where the multiplication is defined using the isomorphism $\omega^{\otimes
    n}_\Delta\cong \Delta$. Then:
\begin{enumerate}
\item $\Lambda$ is a prime    homologically homogeneous ring;
\item $\omega_\Lambda\cong \Lambda$, as $\Lambda$-bimodules.
\end{enumerate}
\end{proposition}

\begin{proof}
(1) By Theorem~\ref{ref-2.3-6}(3,4),  $Z$ is an affine normal domain, 
and $\Delta$ is a tame $Z$-order in its simple artinian ring of fractions $A$.   
  By \cite[Corollary~3.6]{YZ1}, $\omega_\Delta$ is $Z$-central
and so Lemma \ref{new-lemma}(4,5) implies that 
$\omega_\Delta$ is isomorphic to a divisorial fractional ideal $I$. Therefore,
$I^{(n)}=\Lambda a$ for some $a\in L=Z(A)$
 and so   $\Lambda\cong \Delta[I]/(1-ax^n)$. The field of fractions of $\Lambda$ is therefore
\[
B=A\otimes_L L [x]/(1-ax^n)
\]
By \cite[Theorem~2.3]{silver}, $\Div(\Delta)$ is a free abelian group.
 Therefore, if $a=b^m$ for some  $m>1$ and $b\in L$ then 
we would   have $m\mid n$ and $I^{(n/m)}=\Lambda b$ contradicting
the minimality of $n$. If follows that $L[x]/(1-ax^n)$ is a field and thus
$B$ is a central simple algebra. Consequently $\Lambda$ is prime. 

The ring    $\Lambda$ is  strongly graded and hence  $\gldim \Lambda
= \gldim \Delta$ follows   from \cite[Corollary~7.6.18]{MR} together with
the  fact that the categories of $\Delta$-modules and graded $\Lambda$-modules are equivalent.
Thus    $\gldim \Lambda=\gldim \Delta=\GKdim \Delta = \GKdim\Lambda$
by Theorem~\ref{ref-2.3-6}(2) and \cite[Proposition~8.2.9(ii)]{MR}.  
   By Theorem~\ref{ref-2.3-6}(1),  $\Delta$ 
 is CM as a $Z$-module  and hence so is each $\omega_\Delta^{\otimes j}$ 
 and  $\Lambda$. Thus $\Lambda$  is
homologically homogeneous by Lemma~\ref{ref-2.4-7}.

(2) Using the formula $\omega_\Lambda=\Hom_R(\Lambda,\omega_R)$ we
compute  that 
\[
\omega_\Lambda=\omega_\Delta\oplus \omega^{\otimes 2}_\Delta\oplus\cdots\oplus
\omega^{\oplus n-1}_\Delta\oplus \Delta
\]
as $\ZZ/n\ZZ$-graded $\Lambda$-bimodules.  Forgetting the grading gives the result.
\end{proof}

\begin{remarks} \label{ref-3.2-235} 
(1)  Assume that  $Z$ an essentially affine $k$-algebra. 
Following \cite{Br3} or \cite{Gi},   $\Delta$ 
is  called  \emph{Calabi-Yau of  dimension $d$} if   
 $D_\Delta\cong \Delta[d]$  in  $D(\Delta^e)$.    
 (Some authors  also require 
  Calabi-Yau algebras to have finite global dimension; see, for example,  
  \cite[Theorem~3.2(iii)]{IR}.) For a survey on the   Calabi-Yau property in an 
  algebraic context see \cite{Gi}.

 By Proposition~\ref{ref-2.15-18}, an affine  Calabi-Yau algebra of finite global 
 dimension is automatically homologically homogeneous.
Conversely, Proposition~\ref{ref-3.1-23} can be regarded as a reduction to the 
Calabi-Yau case.

(2) Proposition~\ref{ref-3.1-23} can also be regarded as a reduction to the case of
orders unramified in codimension one. In order to explain this, recall that a tame order  $\Delta$  is 
  \emph{unramified in codimension one} if   $\p\Delta_\p = \rad \Delta_\p$ for all
  $p\in \mathfrak{P}_1(Z)$.  Given a tame  Caladi-Yau  order $\Delta$, then  
   Lemma~\ref{new-lemma}(6) implies that  $Z\cong Z(\omega_\Delta)=\omega_Z$ and so  
   Lemma~\ref{new-lemma}(5) implies that  $\p\Delta_\p = \rad \Delta_\p$ for all  
   $p\in \mathfrak{P}_1(Z)$.
 \end{remarks}

Even when $\Delta$ is homologically homogeneous, there is no reason for 
$\omega_\Delta$ to have finite order and so Proposition~\ref{ref-3.1-23}
cannot be applied directly. However, $\omega_\Delta$ has finite order
 locally, which will be sufficient for our applications.
 Before stating the result, we prove some elementary facts.

\begin{lemma}\label{ref-3.2-24}
 If $S$ is a ring with Jacobson radical $\rad(S)$ and $P$
is an invertible $S$-bimodule then $\rad(S) P=P \rad(S)$. 
\end{lemma}

\begin{proof} We claim that the image of composition 
\begin{equation}
\label{ref-3.3-25}  \chi:
P^{-1}\otimes_S \rad(S) \otimes_S P\ \r \
P^{-1}\otimes_S S \otimes_S P\ \cong \   S
\end{equation}
lies in $\rad(S)$. This proves the inclusion $\rad(S) P\subseteq P \rad(S)$. 
To prove the opposite inclusion interchange $P$ and $P^{-1}$.

In order to prove the claim we will  prove that the image of  $\chi$
annihilates all simple $S$-modules.  Let $M$ be a simple $S$-module. 
We must show that the map
\begin{equation}
\label{ref-3.4-26}
P^{-1}\otimes_S \rad(S) \otimes_S P\otimes_S M \r M
\end{equation}
is  zero. 
Tensoring \eqref{ref-3.4-26} on the left by $P$ we obtain the map
\begin{equation}
\label{ref-3.5-27}
\rad(S) \otimes_S P\otimes_S M\r P\otimes_S M
\end{equation}
Since $P\otimes_S-$ is an autoequivalence of $\Mod(S)$,
$P\otimes_S M$ is a simple module  and hence  \eqref{ref-3.5-27}
is indeed the zero map.
\end{proof}

\begin{lemma}   \label{ref-3.6-28}
 Assume   that   $R$ is local and that 
$\Gamma$ is a  tame $R$-order in $A$, with $R=Z(\Gamma)$. 
If  $P$ is an $R$-central invertible  $\Gamma$-bimodule,   
     then there exists an integer $n>0$ such that $P^{\otimes n}\cong
   \Gamma$ as  $\Gamma$-bimodules.
\end{lemma}

\begin{proof}  
Since $P$ is invertible, tensor powers, symbolic powers and ordinary powers 
 all coincide, so we will drop the tensor product sign.

We first prove that $P^{ n}\cong
   \Gamma$ as  left $\Gamma$-modules.
By Lemma~\ref{ref-3.2-24},   $P/\rad(\Gamma)P$ is an invertible
  bimodule over $\Gamma/\rad(\Gamma)$. Since $\Gamma/\rad(\Gamma)$ is
  semi-simple, it is easy to see that there exists  $n>0$ such that
\[
P^{ n}/ \rad(\Gamma)P^{ n}=
 (P/\rad(\Gamma)P)^{ n}\cong \Gamma/\rad(\Gamma)
\]
 as left  $\Gamma/\rad(\Gamma)$-modules.   By Nakayama's Lemma it follows
that $P^{ n}\cong \Gamma$, again as left $\Gamma$-modules.

Let $K$ denote the fraction field of $R$.
 Since $P$ is $R$-central,  $K\otimes_R P$ is an invertible $A$-bimodule. After
choosing an isomorphism $K\otimes_R P\cong A$ we may assume that $P$
is a divisorial fractional $R$-ideal.  
By \cite[Proposition~II.4.20]{LVV}, some power $P^{(e)}$ of $P$ lies in the 
image of $\Div(Z(\Gamma))$ in $\Div (\Gamma)$; that is,
$P^{(e)}\cong  (\Gamma I)^{\ast\ast}$  for some reflexive ideal $I$ of $R$.
By the last paragraph,  we may also assume that $P^{e} \cong  \Gamma$
  as left $\Gamma$ modules. 

   Now let $u=\rk_R \Gamma$. 
   Then  $\big(\bigwedge^u_R(\Gamma I)\big)^{\ast\ast}\cong \bigwedge_R^u P^{e}
\cong \big(\bigwedge^u_R(\Gamma )\big)^{\ast\ast}$ as $R$-modules. 
 On the other hand, $$\big(\bigwedge^u(\Gamma I)\big)^{\ast\ast}=
 \big((\bigwedge^u \Gamma ) I^u\big)^{\ast\ast} 
=\big(\bigwedge^u\Gamma\big)^{\ast\ast} \big(I^u\big)^{\ast\ast}$$  
and so $\textstyle   \big(\bigwedge^u\Gamma\big)^{\ast\ast} \big(I^u\big)^{\ast\ast} =
\big(  \bigwedge^u(\Gamma )\big)^{\ast\ast}   . $    Cancelling
  $\big(\bigwedge_R^u\Gamma\big)^{\ast\ast}$  
  gives $\big(I^u\big)^{\ast\ast}\cong R$ as $R$-modules. Since
    $P^{e}\cong  (\Gamma I)^{\ast\ast}$ as   $\Gamma$-bimodules  we obtain 
  $P^{eu}\cong \Gamma$  as  $\Gamma$-bimodules.
\end{proof}

\begin{corollary} \label{ref-3.7-29} Suppose that $\Delta$ is
  homologically homogeneous $k$-algebra. 
  Then for every maximal ideal $\m$ of $Z$
  there exist $f\in Z\smallsetminus\m$ and $n>0$ with the property that
  $\omega_{\Delta_f}^{\otimes n}\cong \Delta_f$ as  $\Delta_f$-bimodules.
\end{corollary}

\begin{proof} 
By Proposition~\ref{ref-2.15-18},  $\omega_\Delta$ is invertible and, as has been observed 
before, it is $Z$-central. By Lemma~\ref{new-lemma}(3) and Theorem~\ref{ref-2.3-6}(4),
we can therefore 
 apply  Lemma~\ref{ref-3.6-28} to $P=\omega_{\Delta_\mathfrak{m}}$ and conclude that 
$\omega_{\Delta_\m}^{\otimes n}\cong 
\Delta_\m$ as $\Delta_\m$-bimodules. As usual this 
isomorphism may be ``spread out'' on a neighbourhood of $\m$ in $\Spec Z$. 
\end{proof}


\section{The centre of homologically homogeneous rings}\label{ref-4-30}
In this section we  prove Theorem~\ref{ref-1.1-1} from the introduction. We start with
 two preparatory lemmas, the first of which gives a useful  algebraic criterion for a 
 ring to have rational singularities. 

\begin{lemma}\label{ref-4.3-32}
Let $Z$ be an affine    normal CM $k$-domain with field of fractions $K$. 
Then $Z$ has rational singularities if and only if, for all regular 
 affine $k$-algebras $S$ satisfying  $Z\subseteq S\subset K$, 
we have $\omega_Z\subseteq \omega_S$ inside $\omega_K$.
\end{lemma}

\begin{proof} Let $X=\Spec Z$.  By Remark~\ref{canonical-defn},
$\omega_X$ in the sense of 
\cite{KK, KM} is equal to $\omega_Z$ in the sense of this paper and so, 
by Lemma~\ref{new-lemma}(4), $\omega_X$ is reflexive.
 According to  \cite[p.~50]{KK} or  \cite[Theorem~5.10]{KM}, $X$
  has rational singularities if and only if for one (or for all)
  resolution(s) of singularities $f:Y\r X$ we have
  $f_\ast\omega_Y=\omega_X$ inside $\omega_K$. Since $\omega_X$ and $\omega_Y$ are
reflexive this is equivalent to $\omega_X\subseteq f_\ast\omega_Y$
    and  the latter condition  is equivalent to $(f^\ast \omega_X)^{\ast\ast}\subseteq
  \omega_Y$.  This  can be checked locally on $Y$. 

  So assume that  $\omega_Z\subseteq \omega_S$  for all affine    regular $k$-algebras $S$
satisfying  $Z\subseteq S\subset K$. 
Pick $Y$ by the last paragraph and an open affine subset $U\subset Y$.
 Then   $\omega_Z\subseteq \omega_S$  for $S=\mathcal{O}(U)$  and hence 
$(S\otimes_Z \omega_Z)^{\ast\ast} \subseteq \omega_S$.  
Globalizing gives $(f^\ast \omega_X)^{\ast\ast}\subseteq
  \omega_Y$ and so  $Z$ has rational singularities. 

Conversely assume that $Z$ has rational singularities and let
$Z\subseteq S$ be as in the statement of the lemma. Put $U=\Spec S$. We may
compactify the map $g:U\r X$ to a projective map
$\bar{g}:Y'\r X$. A priori $Y'$ will not be smooth but we can resolve it
further without touching $U$ (see \cite[Theorem~0.2]{KM}) to arrive at a
resolution of singularities $f:Y\r X$. The fact that $(f^\ast \omega_X)^{\ast\ast}
\subseteq \omega_Y$ yields $(S\otimes_Z
\omega_Z)^{\ast\ast}\subseteq \omega_S$  after restricting to $U$. Thus
 $\omega_Z\subseteq \omega_S$. 
\end{proof}

\begin{lemma}\label{ref-4.1-312}  Let $\Lambda_1$ and $\Lambda_2$ be  affine $k$-algebras of 
finite global dimension that satisfy a polynomial identity. Then 
$\Lambda_1\otimes_k\Lambda_2$ has finite global dimension. 
\end{lemma}

\begin{proof} By the Nullstellensatz \cite[Theorem~13.10.3]{MR}, every primitive 
 factor ring of $\Lambda_i$ is isomorphic to a full matrix ring over $k$. Hence every primitive 
  factor  ring $\Gamma$ of $\Lambda=\Lambda_1\otimes_k\Lambda_2$ decomposes 
as $\Gamma=\Gamma_1\otimes_k \Gamma_2$ for primitive factor rings  $\Gamma_i$ of $\Lambda_i$.
Thus  any simple $\Lambda$-module $M$ can be written as $M=M_1\otimes_kM_2$,
where each $M_i$ is a simple $\Lambda_i$-module. Now use 
\cite[Proposition~IX.2.6]{CE}.
\end{proof}

\begin{theorem} \label{ref-4.1-31}  
If  $\Delta$ is a homologically homogeneous $k$-algebra, then
  $Z=Z(\Delta)$ has rational singularities. 
\end{theorem}

\begin{proof} 
 It is enough to prove the result locally, so by 
 Corollary~\ref{ref-3.7-29} we can replace $\Delta$ by some $\Delta_f$ and assume that 
  $\omega_{\Delta}^{\otimes n}\cong \Delta$ as
  $\Delta$-bimodules.  By  Proposition~\ref{ref-3.1-23}, the algebra
  $\Lambda= \Delta\oplus \omega_\Delta\oplus \omega^{\otimes 2}_\Delta\oplus\cdots\oplus
\omega^{\otimes n-1}_\Delta$ satisfies $\omega_\Lambda\cong \Lambda$ as $\Lambda$-bimodules.
Then $\Lambda$ and hence $Z(\Lambda)$ are $\mathbb{Z}/n\mathbb{Z}$-graded.
 Moreover,  as $ \omega_\Delta$ is $Z$-central, clearly $Z$ commutes with 
each $\omega_\Delta^{\otimes j}$ and so $Z\subseteq Z(\Lambda)_0$. Since
the  other inclusion is trivial,  
$Z =  Z(\Lambda)_0$ and $Z$  is a module-theoretic summand of $  Z(\Lambda)$. 
Since a direct summand of a  ring with
  rational singularities has rational singularities \cite{boutot} we
  may therefore   replace $\Delta$ by $\Lambda$ and assume that  
   that $\omega_\Delta\cong\Delta$ as bimodules. 
By Proposition~\ref{ref-3.1-23}(1) $\Delta$ remains homologically homogeneous. 

  We will use Lemma~\ref{ref-4.3-32}, so 
fix a ring  $Z\subseteq S\subset K$ as in the lemma
 and let $\Gamma$ be a maximal, and therefore tame  $S$-order containing
$S\Delta$ inside the simple artinian ring of fractions $A$ of $\Delta$.
  Our discussion in Section~\ref{ref-2-2} on dualizing complexes also applies to $\Gamma$,  so
$\omega_\Gamma=\Hom_S(\Gamma,\omega_S)$ and $\omega_S=Z(\omega_\Gamma)$
in the notation developed there.
   We will show that $\omega_\Delta\subseteq\omega_\Gamma$
inside $\omega_A$.  Since $S\subset K$, this will yield $Z(\omega_\Delta)\subseteq
Z(\omega_\Gamma)$ as subgroups of $Z(\omega_A)$ 
and so Lemma~\ref{new-lemma}(6)  will imply that  $\omega_Z\subseteq \omega_S$, as required.

\medskip

In order to prove that $\omega_\Delta\subseteq\omega_\Gamma$ we may as
well prove that $\omega_\Delta\Gamma\subseteq \omega_\Gamma$.  The
bimodule isomorphism $\omega_\Delta\cong \Delta$ means that
$\omega_\Delta=c\Delta$ for some central element $c\in \omega_A$. From this
we deduce that $\omega_\Delta\Gamma=\Gamma\omega_\Delta$ is an
invertible $\Gamma$-bimodule with inverse $\omega^{(-1)}_\Delta \Gamma=c^{-1}\Gamma$.   
By Lemma~\ref{new-lemma}(4), $\omega_\Gamma$ is reflexive and so 
 it    suffices to prove
that $\omega^{(-1)}_\Gamma\subseteq\omega^{(-1)}_\Delta \Gamma$ inside
$\omega_A^{(-1)}$.

We claim that 
\begin{equation}\label{Gamma-equ}
\Gamma\Lotimes_{_{\scriptstyle\Delta}}\hskip -2pt
 \RHom_{\Delta^e}(M,\Delta^e) \Lotimes_{_{\scriptstyle\Delta}}\hskip -2pt \Gamma=
\RHom_{\Gamma^{e}}(\Gamma\Lotimes_{_{\scriptstyle\Delta}}\hskip -2pt
 M \Lotimes_{_{\scriptstyle\Delta}}\hskip -2pt\Gamma,\Gamma^e)
\end{equation}
for any object $M$ in $D^b(\Delta^e)$ with finitely generated cohomology.  
To prove this recall that, by Lemma~\ref{ref-4.1-312},   $\gldim \Delta^e<\infty$.
 Thus we can  replace $M$ by a finite projective  resolution of $\Delta^e$-modules and
it then  suffices prove the claim for  $M=\Delta^e$. This case is obvious.

Applying \eqref{Gamma-equ}  with $M=\Delta$ and using  the formula 
$\omega_\Delta^{(-1)}=\RHom_{\Delta^e}(\Delta,\Delta^e)[d]$ from 
\eqref{ref-2.19-22} we obtain
\[
\Gamma\Lotimes_{_{\scriptstyle\Delta}}\hskip -2pt \omega_\Delta^{(-1)}
\Lotimes_{_{\scriptstyle\Delta}}\hskip -2pt \Gamma=
\RHom_{\Gamma^{e}}(\Gamma\Lotimes_{_{\scriptstyle\Delta}}\hskip -2pt\Gamma,\Gamma^e)[d].
\]
Using the fact that the derived tensor product maps to the ordinary tensor
this induces a composed map
\[
\RHom_{\Gamma^{e}}(\Gamma,\Gamma^e)[d] \ \r \  
\RHom_{\Gamma^{e}}(\Gamma\Lotimes_{_{\scriptstyle\Delta}}\hskip -2pt\Gamma,\Gamma^e)[d]
\ = \  \Gamma\Lotimes_{_{\scriptstyle\Delta}} \omega^{(-1)}_\Delta
\Lotimes_{_{\scriptstyle\Delta}}\hskip -2pt \Gamma \  \r  \
\omega^{(-1)}_\Delta\Gamma.
\]
Taking cohomology in degree zero and then biduals  gives a map
\[
\Ext^d_{\Gamma^e}(\Gamma,\Gamma^e)^{\ast\ast}\r ( \omega^{(-1)}_\Delta\Gamma)^{\ast\ast} = 
 \omega^{(-1)}_\Delta\Gamma.
\]
Using \eqref{ref-2.18-21} this induces a map
\begin{equation}
\label{ref-4.4-33}
\omega_\Gamma^{(-1)}\r \omega_\Delta^{(-1)}\Gamma.
\end{equation}
Now we could have done these computations after tensoring with the field of fractions $K$ of $Z$.
Since $K=KZ=KS$ and $K\Delta=K\Gamma=A$, all morphisms would then have been
(canonically) the identity.  From this we deduce that \eqref{ref-4.4-33} is
an inclusion which takes place inside $\omega_A^{(-1)}$. This means we are done. 
\end{proof}

\begin{remarks}\label{ref-4.2-31} (1)  
Suppose that $\Delta$ is an affine Calabi-Yau $k$-algebra of finite global dimension. 
Then Theorem~\ref{ref-4.1-31} and  Remark~\ref{ref-3.2-235}(1)
combine to prove that   $Z$  has rational singularities.

(2)   Homologically homogeneous rings were defined in \cite{BH} for orders in 
semisimple  rather than simple artinian  rings. However, by   \cite[Theorem~5.3]{BH}, 
 these more general  algebras are  direct sums of
  prime homologically homogeneous rings and so the more general case also follows from this 
  theorem.  Similarly, one can weaken the hypothesis that $\Delta$ be  finitely generated as a module 
over its centre to the assumption that it be an affine algebra
satisfying a polynomial identity
 since, by \cite[Theorem~5.6(iv)]{SZ}, this 
already forces $\Delta$ to be a  finitely generated  $Z$-module.  
\end{remarks}


\section{Examples}\label{ref-5-34}
Here we give two examples to illustrate the earlier results. 
The first shows that \cite[Lemma~4.2]{VdB32} cannot be improved while 
the second one  shows that Theorem~\ref{ref-1.1-1} can fail in finite 
characteristic.

In addition to our standing hypotheses, 
suppose that $R$ is  an affine Gorenstein $k$-algebra and that
 $\Delta=\End_R(M)$ for some finitely generated reflexive $R$-module $M$. Then
it follows from  \cite[Lemma~4.2]{VdB32} 
that $\Delta$ is homologically homogeneous if and only if 
  $\gldim \Delta<\infty$  and $\Delta$ is a CM $R$-module. This is 
useful for   the theory of  noncommutative crepant resolutions, so it 
would be useful if we could weaken the hypotheses in this result.  In our 
first example, we show that the Gorenstein condition is necessary.

\medskip
Here is the example. Let $T$ be a one-dimensional torus acting on the
generators of the polynomial ring $S=k[x_1,x_2,x_3,x_4,x_5]$ with
weights $1,1,1,-1,-1$ and let $R=S^T$.  We may also view $R$  as the
coordinate ring of the variety of $2\times 3$-matrices of rank $\le 1$.

The $T$-weights give a grading   $S=\bigoplus_{\ell=-\infty}^\infty S_\ell$ 
with $S_0=R$. According to  the proof of \cite[Lemma~8.8]{VdB32} the $S_i$ 
are isomorphic to reflexive ideals of $R$ with $S_{a+b}=(S_aS_b)^{\ast\ast}$
for all $a,b\in \mathbb{Z}$.
Furthermore it is easy to see that  $S_i$ is not  a projective 
$R$-module when $i\not= 0$. 

It follows from \cite[Lemma~8.1]{VdB32} that $S_{-2},S_{-1}, R$ and $S_1$
are CM $R$-modules while $R$ is certainly normal. 
 It follows from  \cite[Lemma~8.2 and Theorem~8.6]{VdB32}
 that 
\[
\Delta=\End_{R}(R \oplus S_1)=
\begin{pmatrix}
R  \hfill & S_1 \hfil \\
S_{-1}  \hfil& R \hfill
\end{pmatrix}.
\]
has finite global dimension and hence is a tame order over its centre
$R$. By \cite[Korollar 2]{Knop1}, the dualizing module $\omega_R$ 
is  isomorphic to $S_{-1}$ (where $-1$ represents minus the sum of the
weights of the generators of $S$)   from which we deduce that
\[
\omega_\Delta = \Hom_R(\Delta,\omega_R) \cong 
\begin{pmatrix}
S_{-1}  \hfil& R \hfill \\
S_{-2}  \hfil & S_{-1} \hfil
\end{pmatrix}.
\]
Both $\Delta$ and $\omega_\Delta$ are graded for the standard grading on $R$.
For this choice of grading, $\Delta$ is graded semi-local and $\omega_\Delta$ is (as left
or right module) not a direct sum of indecomposable graded $\Delta$-projectives.
Consequently, $\omega_\Delta$ is not projective.

By Proposition~\ref{ref-2.15-18},  $\Delta$ is therefore not homologically 
homogeneous. 

\begin{remarks} (1) By the proof of \cite[Proposition~A.2]{DNVdB} it follows that 
$\omega_\Delta$ defines an element of the \emph{derived} Picard group of 
$\Delta$.

(2) The methods of \cite{VdB34} allow one to treat this example
in the context of determinantal varieties. It follows from the results given there
that one of the  simple  graded $\Delta$-modules  has projective dimension $4$
and the other has projective dimension $5$.
\end{remarks}

The example leads naturally to the following question.

\begin{question}  
 Assume that $Z=Z(\Delta)$ is an affine normal $k$-domain and that
 $\Delta$ is a finitely generated CM  $Z$-module
  with finite global  dimension.  Then, does $Z$ have rational singularities?
\end{question}

We now turn to an example in finite characteristic of a homologically
homogeneous ring whose centre is CM but which does not have
rational singularities in any reasonable sense.

  Assume that $F$ is a field of characteristic $2$ and let   $C=F[u,v,x,y]/(p,q)$
where
$$p=x+u^2+x^2u\qquad\mathrm{and}\qquad
q=y+v^2+y^2v.$$
As the Jacobian matrix of $p,q$ with respect to $x,y$ is invertible,
  $C/F[x,y]$ is \'etale  and hence $C$ is regular. Consider the
action of $G=\mathbb{Z}/(2)=\{1,\sigma\}$ on $C$ by $\sigma(u)=u+x^2$, $\sigma(v)=v+y^2$,
$\sigma(x)=x$, and $\sigma(y)=y$. 
 Clearly $B=C^G$ is an affine normal domain of Krull
dimension two and hence it is CM.  

Resolutions of singularities are known exist for surfaces in all
characteristics and there is a corresponding satisfactory theory of
rational singularities. We will  show that $B$ does not have
rational singularities.   Let $\m=(u,v)\subset C$ and notice that $\m = (u,v,x,y)$
is maximal; thus $\widehat{C}_\m = F[[u,v]]$.    It suffices to
prove that $\widehat{B}_\mathfrak{n}$, for $\mathfrak{n}=B\cap \m$, does not
have rational singularities.  Since $ uu^\sigma= u^2+ux^2=x\in
\widehat{C}^G_\m=\widehat{B}_\mathfrak{n}$ and $vv^\sigma=y\in
\widehat{B}_\mathfrak{n}$, our notation conforms with that of
\cite[Theorem]{Ar}.  Now the fact that $u^2+x^2u+x=0=v^2+y^2v+y$ means
that $\widehat{B}_\mathfrak{n}$ does not have rational singularities
by the observation from \cite[p.~64]{Ar}.

Finally, let  $\Lambda=C[x;\sigma]$ be the twisted polynomial ring; thus $xc=c^\sigma
x$ for all $c\in C$.  By the Nullstellensatz, every simple
$\Lambda$-module is finite dimensional and so, by \cite[Theorem~7.9.16]{MR},
$\Lambda$ is homologically homogeneous of dimension $3$.  As
$\sigma^2=1$, the element $x^2$ is central. It follows routinely that
$Z(\Lambda)=B[x^2]$. Thus, $Z(\Lambda)$ also does not have rational  singularities.

The basic reason why such counterexamples exist in bad characteristic is that 
a fixed ring $S^G$ need not be a summand of the ring $S$.
The  example \cite[Example~7.3]{BHM} of a homologically homogeneous ring 
with a non-CM centre  occurs for a similar reason. So, it is natural to ask:

\begin{question} Suppose that $\Lambda$ is a homologically homogeneous ring 
whose centre $Z(\Lambda)$   is an 
affine  $F$-algebra for field $F$ of characteristic $p>0$. 
If $Z(\Lambda)$ is a $Z(\Lambda)$-module  summand 
of $\Lambda$, then does $Z(\Lambda)$ have rational singularities?
 \end{question}

It was conjectured in \cite{VdB-crepant} and proved in \cite[Theorem~6.6.3]{VdB32}
that a $3$-dimensional $k$-variety  with terminal singularities 
has a noncommutative crepant resolution if and only if it has a commutative one
(see also  \cite[Corollary~8.8]{IR}).
We end by noting that this is not true in higher dimensions.
One way to produce counterexamples is with the fixed ring 
$R=\mathbb{C}[V]^G$
of a finite group $G\subset SL(V)$, where $V=\mathbb{C}^n$.
In this case,  the  twisted group ring $\mathbb{C}[V]\ast G\cong
  \End_{R}(\mathbb{C}[V])$ is a noncommutative crepant resolution of 
$R$ \cite[Example~1.1]{VdB32}, but it is well-known that such a ring $R$ 
need not have a commutative crepant resolution (see, for example, 
\cite[Theorem~1.7]{Kn}).

\def\cprime{$'$} \def\cprime{$'$}
\ifx\undefined\bysame
\newcommand{\bysame}{\leavevmode\hbox to3em{\hrulefill}\,}
\fi

\end{document}